\documentclass[11pt]{article}
\usepackage{amssymb}
\usepackage{theorem}

\newtheorem{theorem}{Theorem}
\newtheorem{prop}[theorem]{Proposition}
\newtheorem{lemma}[theorem]{Lemma}

\newtheorem{corollary}[theorem]{Corollary}
\newtheorem{remark}[theorem]{Remark}

\theorembodyfont{\rmfamily}
\newtheorem{example}{Example}
\newenvironment{proof}{\noindent{\bf Proof. }}{\hfill$\square$\medskip}

\def\R{{\mathbb R}}
\def\one{{\mathbf 1}}

\def\eps{\varepsilon}

\def\FF{{\cal F}}
\def\AA{{\cal A}}
\def\BB{{\cal B}}

\def\LL{{\cal L}}

\def\CC{{\cal C}}

\def\E{{\sf E}}
\def\T{^{\sf T}}

\begin{document}
\begin{center}

{\Large\bf
Determining the Genus of a Map\\[6pt]by Local Observation\\[6pt]of a Simple
Random Process
}\\[12mm]
{\sc Itai Benjamini} and {\sc L\'{a}szl\'o Lov\'{a}sz}\\
Microsoft Research \\
One Microsoft Way, Redmond, WA  98052\\
e-mail: itai@wisdom.weizmann.ac.il lovasz@microsoft.com \\[20mm]

\end{center}

\begin{abstract}
Given a graph embedded in an orientable  surface, a process
consisting  of random excitations and random
node and face balancing is constructed and analyzed. It is shown
that given a priori bounds $\overline{g}$ on the genus and
$\overline{n}$ on the number of nodes, one can determine the genus of
the surface from local observations of the process restricted to any
connected subgraph which cannot be separated from the rest of the
graph by fewer than $16\overline{g}$ nodes. The observation time and
the computation time are polynomial in $\overline{n}^{\overline{g}}$.

The process constructs slightly perturbed random ``discrete analytic
functions'' on the surface, and the key fact in the analysis is that
such a function cannot vanish on a large piece of the surface.
\end{abstract}

\section{Introduction}
At least since Polya (1921) proved that ``a drunk man will
return home while a drunk bird might lose its way forever'', it is
known that geometric properties of the underling space can manifest
themselves in the behavior of a random process taking place on the
space. Moreover, in recent years random processes were used to
retrieve information about the underlying space. E.g. sampling,
volume estimates of convex bodies or scenery reconstruction along
random walk paths (see for instance \cite{LW}, \cite{Bo} and
\cite{Ke}). Rather then letting a random walker wander around and
gather information, in this note we would like to study a problem in
which a random process is observed locally in a fixed bounded
neighborhood and still non-trivial global observations can be
distilled from these observations.

In addition to the algorithmic motivation, a reasonable question
is as follows. Consider a stationary spin system on a graph, such as
Glauber dynamics for the Ising model on a graph. What properties of
the graph can be inferred from properties of the process (for
instance, are there interesting relations between the spectrum of
the Ising dynamics and the graph spectrum)?  A harder challenge is
to infer non-trivial global properties of the underling graph from
local observations of the process.

Facing the harder challenge one might devise first custom made
variants of the standard spin systems, which can be analyzed and
provide ways to compute global invariants using ``physical" systems.

Indeed, below we will study a reasonably natural and simple process,
called {\it noisy circulator}, with local operations, living on the
edges of a graph, embedded in an orientable surface. The noisy
circulator consists of adding mass $1$ to randomly chosen edges with
a slow rate and balancing the flow into vertices or around faces at
random and with a faster rate (details below). Somewhat surprisingly,
it will be shown how to extract, with high probability and in polynomial time
in the size of the graph, the genus of the surface by observing the
restriction of the process to a bounded set of edges, depending only
on an a priori bound on the genus.

Although from a pure algorithmic view point this first construction
might be useful and might have some advantages, the point is not to
devise the optimal ad-hoc distributed algorithm for finding the
genus, under certain restrictions, but to show how  locally observing
a simple physical process already does that, and to present this
computation scheme. It is of interest then to construct other
examples of similar flavor.

The structure of the proof is twofold, a topological theorem and a
statistical element.

In the topological part, we study {\it discrete holomorphic
1-forms} which we refer to here as {\it smooth circulations}:
circulations that are also circulations on the dual map. These can be
considered as discrete analogues of analytic functions. Every
homology class of circulations contains exactly one smooth
circulation, so the dimension of their space can be used to find the
genus of the surface. A key result (Theorem~\ref{NONDEG}), which is
of independent interest, asserts that every connected piece of the
vanishing set of a smooth circulations can be separated from the rest
of the graph by a small number of nodes.

In the statistical part, we determine the dimension of the space of
smooth circulations from observations which can be considered as
samples of smooth circulations restricted to a bounded set of edges
with additional noise.

The next subsection contains some necessary definitions and basic
propositions; in subsection~\ref{RCIRC} a description of the process
and a formulation of the main theorem are given, together with an
outline of the proof. The rest of the paper contains the proof of
the main theorem.  In the final section we end with some further
comments and problems.

\subsection{Circulations and homology}

Let $S$ be a closed compact surface, and consider a {\it map} on
$S$, i.e., a graph $G=(V,E)$ embedded in $S$ so that each face is a
disc. We can describe the map as a triple $G=(V,E,\FF)$, where $V$ is
the set of nodes, $E$ is the set of edges, and $\FF$ is the set of
faces of $G$. We fix a reference orientation of $G$; then each edge
$e\in E$ has a {\it tail} $t(e)\in V$, a {\it head} $h(e)\in V$, a
{\it right shore} $r(e)\in \FF$, and a {\it left shore} $l(e)\in \FF$.

The embedding of $G$ defines a {\it dual map} $G^*$. Combinatorially,
we can think of $G^*$ as the triple $(\FF,E,V)$, where the meaning of
``node'' and ``face'', ``head'' and ``right shore'', and ``tail''
and ``left shore'' is interchanged. (Taking the dual of the dual
will give the original map with every edge reversed; this should not
concern us in this paper.)

For each edge $e$, let $\chi_e\in\R^E$ be the unit vector that is $1$
on $e$ and $0$ elsewhere; we define $\chi_v\in\R^V$ for $v\in V$ and
$\chi_F\in\R^\FF$ for $F\in \FF$ analogously.

For each node $v$, let $\delta v\in\R^E$ denote the coboundary of
$v$:
\[
(\delta v)_e=\cases{1,& if $h(e)=v$,\cr
                  -1,& if $t(e)=v$,\cr
                  0,& otherwise.\cr}
\]
Thus $|\delta v|^2=d_v$ is the degree of $v$. A vector $\phi\in\R^E$
is a {\it circulation} if
\[
\phi\cdot \delta v = \sum_{e:~h(e)=v} \phi(e)-\sum_{e:~t(e)=v}
\phi(e)=0.
\]

For every face $F\in \FF$, we denote by $\partial F\in\R^E$ the
boundary of $F$:
\[
(\partial F)_e=\cases{1,& if $r(e)=F$, \cr
                  -1,& if $l(e)=F$, \cr
                  0,& otherwise.\cr}
\]
Then $d_F=|\partial F|^2$ is the length of the cycle bounding $F$.

Each vector $\partial F$ is a circulation; circulations that are
linear combinations of vectors $\partial F$ are called {\it
null-homologous}. Two circulations $\phi$ and $\phi'$ are called {\it
homologous} if $\phi-\phi'$ is null-homologous.

Let $\phi$ be a circulation on $G$. We say that $\phi$ is {\it
smooth} if for every face $F\in \FF$, we have
\[
\phi\cdot \partial F=0.
\]
This is equivalent to saying that $\phi$ is a circulation on the dual
map $G^*$.

\begin{remark}
Smooth circulations are closely related to {\it discrete analytic
functions} and are essentially the same as {\it discrete holomorphic
1-forms}. These functions were introduced for the case of the square
grid a long time ago \cite{Fe,Du}. For the case of a general planar
graph, the notion is implicit in \cite{Br}. For a detailed treatment
see \cite{Me}.

To explain the connection, let $\phi$ be a smooth circulation on a
graph $G$ embedded in a surface. Consider a planar piece of the
surface. Then on the set ${\cal F}'$ of faces contained in this
planar piece, we have a function $\sigma:~{\cal F}'\to \R$ such that
$\partial\sigma=\phi$, i.e., $\phi(e)=\sigma(r(e))-\sigma(l(e))$ for
every edge $e$. Similarly, we have a function $\pi:~V'\to\R$ (where
$V'$ is the set of nodes in this planar piece), such that
$\delta\pi=\phi$, i.e., $\phi(e)=\pi(t(e))-\pi(h(e))$ for every edge
$e$. We can think of $\pi$ and $\sigma$ as the real and imaginary
parts of a (discrete) analytic function. The relation
$\delta\pi=\rho\phi$ is then a discrete analogue of the
Cauchy--Riemann equations.
\end{remark}

Thus we have the two linear orthogonal subspaces $\AA\subseteq\R^E$
generated by the vectors $\delta v$ ($v\in V$) and $\BB\subseteq
\R^E$ generated by the vectors $\partial F$ ($F\in \FF$). Vectors in
$\BB$ are 0-homologous circulations. The orthogonal complement
$\AA^\perp$ is the space of all circulations, and $\BB^\perp$ is the
space of circulations on the dual graph. The intersection
$\CC=\AA^\perp\cap\BB^\perp$, the space of smooth circulations. So
$\R^E=\AA\oplus\BB\oplus\CC$. From this picture we conclude the
following.

\begin{prop}\label{SMOOTHHOM}
Every circulation is homologous to a unique smooth circulation.
\end{prop}

It also follows that $\CC$ is isomorphic to the first homology group
of $S$ (over the reals), and hence we get the following:

\begin{prop}\label{SMOOTHDIM}
The dimension of the space $\CC$ of smooth circulations is $2g$.
\end{prop}

\subsection{Main result: randomized circulations}\label{RCIRC}

\subsubsection{The Noisy Circulator.}

We consider the following process on $G$. Let $p>0$ be fixed. We
start with the vector $x=0\in\R^E$. At each further step, the
following two operations are carried out on the current vector
$x\in\R^E$:

\smallskip

(a) [Node balancing.] We choose a random node $v$, and subtract from
$x$ the vector $(x\T\delta_v/d_v)\delta_v$.

\smallskip

(b) [Face balancing.] We choose a random face $F$, and subtract from
$x$ the vector $(x\T\partial_F/d_F)\partial_F$.

\smallskip

In addition, with some given probability $p>0$, we do the following:

\smallskip

(c) [Excitation.] We choose a random edge $e$, and add $1$ to $x_e$.

\smallskip

Immediately after a node balancing step, the node $v$ just balanced
satisfies the flow condition; a subsequent other node balancing, or a
face balancing may destroy this. Smooth circulations are invariant
under node and face balancing, and we'll see that under repeated
application of (a) and (b), any vector converges to a smooth
circulation.

\subsubsection{Observing the process.}

Let $U\subseteq V(G)$ induce a connected subgraph of $G$, and assume
that every cycle that separates $S$ into two parts $S_1$ and $S_2$
so that $S_1$ contains $U$ and $S_2$ is not a disc, has length at
least $16g$. We show that if $p$ is small enough, then observing the
process on the edges incident with $U$ long enough, we can determine
the genus of the surface.

To be precise, we do the following. Let $E_0$ be the set of edges
incident with $U$. Let $x(t)\in \R^E$ be the vector after $t$ steps,
and let $y(t)$ be the restriction of $x(t)\in\R^{E_0}$ to the edges
in $E_0$. So we can observe the sequence random vectors $y(0), y(1),
\dots$.

Let $x'(t)$ be the projection of $x(t)$ onto $\CC$, and let $y'(t)$
be the restriction of $x'(t)$ to the edges in $E_0$. Because of the
random steps (c), after a sufficiently long time, the vectors
$x'(0),\dots,x'(t-1)$ will span the space $\CC$. So the rank of this
set of vectors gives us the dimension of $\CC$, and so by Proposition
\ref{SMOOTHDIM}, it gives us the genus of the surface. The
restriction to $E_0$ is one-to-one on $\CC$ (Theorem \ref{NONDEG}),
and so this rank is the same as the rank of the set $\{y'(0),
y'(1),\dots, y'(t)\}$.

Unfortunately, we cannot observe the vectors $y'(t)$, only the
vectors $y(t)$. But (at least if $p$ is small) we expect $y(t)$ to be
close to $y'(t)$. Indeed, the ``errors'' $x''(t)=x(t)-x'(t)$ are in
$\AA\oplus\BB$, and therefore it is not hard to show that they tend
to $0$ exponentially fast (see lemma \ref{CONVERGE} below), at least
as long as no ``excitation'' step (c) occurs. Therefore the
restrictions $y''(t)=y(t)-y'(t)$ also tend to 0 exponentially. The
speed of convergence depends on the eigenvalue gap of the transition
matrix of the (undirected) random walk on $G$ and $G^*$.

So we are lead to the following standard statistical problem: there
is a sequence $y'(0), y'(1), \dots$ of vectors in $\R^k$, which span
a linear subspace $L$. We observe the sequence $y(t)=y'(t)+y''(t)$,
where the ``error'' $y''(t)$ is small on the average. We want to find
the dimension $k$.

The random vectors $x(t)$ or $x'(t)$ are not independent; but if we
take the differences, i.e., we look at the vectors $x'(t+1)-x'(t)$,
then these are independent (in the probabilistic sense). Indeed,
node balancing and face balancing don't change $x'$; so if no
excitation occurred in step $t+1$, then $x'(t+1)-x'(t)=0$, and new
flow was created on edge $e$, then it depends only on $e$. Hence the
vectors $y'(t+1)-y'(t)$ are also mutually independent.

Since $y'(t+1)-y'(t)=0$ with probability $1-p$, it makes sense to
aggregate $N=1/p$ of these terms to one. So we consider the vectors
$z(t)=y(Nt)-y(N(t-1))$ and $z'(t)=y'(Nt)-y'(N(t-1))$. Then the
vectors $z'(t)$ are mutually independent samples from some
distribution on $L$. (The errors $z''(t)=z(t)-z'(t)$ may be
dependent.)

A further difficulty is that if an excitation step occurs close to
the end of an aggregated interval $[N(t-1)+1,Nt]$, then the error
$z''(t)$ can be larger than the main term $z'(t)$. This happens with
small but not negligible probability. We handle this by randomly
selecting just a fraction of these intervals, so that the
probability of any of these bad large errors occurring is small.

To formalize, we propose the following algorithm to recover the
dimension of $\CC$.

\medskip

\noindent{\bf Genus estimate.} We assume that we are given upper
bounds $\overline{n}\ge n+m+f$ and $\overline{g}\ge g$. Let $m_0$ be
the number of edges incident with $U$. Set
\[
T'=6(\overline{n}+\overline{g}), \qquad
\eps=\overline{n}^{-m_0\overline{n}}, \qquad T=4T'^2.
\]
Construct a sequence of integers $t_1,t_2,\dots,\in[0,T-1]$ as
follows. If we have $t_1,t_2,\dots,t_k$, then compute the linear hull
$\LL(k)$ of $z(t_1), z(t_2),\dots, z(t_k)$. Let $H(k)$ be the set of
integers $t\in[0,T-1]$ for which the vector $z(t)$ is farther from
$\LL(k)$ than $\eps$. If $|H(k)|< T'$, then return $k/2$ as your
guess for $g$. Else, choose a number $t_{k+1}\in H(k)$ randomly and
uniformly. If $k$ becomes larger that $2\overline{g}$, declare the
procedure a failure and stop.

\medskip

The main result of this paper is the following.

\begin{theorem}\label{OBSERVE}
The Genus Estimate Algorithm returns the correct genus with
probability at least $2/3$.
\end{theorem}

If you find that a success probability of $2/3$ is not reassuring
enough, independent repetition of the observation can boost this
arbitrarily close to 1.

\section{Properties of smooth circulations}\label{SMOOTHPROP}

\subsection{Harmonic functions}

It is an easy well-known fact that if $G=(V,E)$ is a connected graph
(the orientation is not needed right now), then for any two nodes
$a,b\in V$ there is a vector $\pi\in\R^V$ such that for every node
$v$,
\begin{equation}\label{HARMONIC}
\sum_{u:~uv\in E} \pi_u - d_v\pi_v
 =\cases{1, & if $v=b$, \cr
        -1, & if $v=a$, \cr
        0, & otherwise.\cr}
\end{equation}
This last expression is equivalent to saying that $f_e=
\pi_{h(e)}-\pi_{t(e)}$ is a flow. We denote this vector $\pi$ by
$\pi_{a,b}$ if we want to express that it depends on $a$ and $b$; if
there is an edge $e$ with $h(e)=b$ and $t(e)=a$, then we also denote
$\pi_{a,b}$ by $\pi_e$.

The vector $\pi$ is not unique; we can add the same scalar to each
entry. For our purposes, it will be convenient to choose it so that
\begin{equation}\label{STANDARD}
\sum_u \pi_u=0.
\end{equation}

There are many interpretations of these functions; for example, let
the graph represent an electrical network, with the edges having unit
resistance. Send a unit electric current from $b$ to $a$. Then the
potential of $u$ is $\pi_u$. The function $\pi_{a,b}$ is often called
a {\it harmonic function with poles $a$ and $b$}.

In terms of the Laplacian $L$ of the graph, the equations
(\ref{HARMONIC}) can be written as
\[
L\pi = \chi_b-\chi_a.
\]
The matrix $L$ is not quite invertible, but it has a one-dimensional
nullspace spanned by the vector $\one=(1,\dots,1)\T$, and so it
determines $\pi$ up to adding the same scalar to every entry. We
assumed in (\ref{STANDARD}) that $\one\T \pi=0$. If $J\in\R^{V\times
V}$ denotes the all-1 matrix, then
\[
(L+J)\pi = L\pi = \chi_b-\chi_a,
\]
and so we can express $\pi$ as
\begin{equation}\label{PI}
\pi=(L+J)^{-1}(\chi_b-\chi_a).
\end{equation}

We can use harmonic functions to give a more explicit description of
smooth circulations in a special case. For any edge $e$ of $G$, let
$\eta_e$ be the projection of $\chi_e$ onto $\CC$.

\begin{lemma}\label{SMOOTHFORM}
Let $a,b\in E$ be two edges of $G$. Then
\[
(\eta_a)_b=\cases{(\pi_b)_{h(a)}-(\pi_b)_{t(a)} +
                     (\pi^*_b)_{r(a)}-(\pi^*_b)_{l(a)} +1, &
                        if $a=b$,\cr
                    (\pi_b)_{h(a)}-(\pi_b)_{t(a)} +
                     (\pi^*_b)_{r(a)}-(\pi^*_b)_{l(a)},&
                        otherwise.\cr}
\]
\end{lemma}

\begin{proof}
Let $x_1$, $x_2$ and $x_3$ be the projections of $\chi_b$ on the
linear subspaces $\AA$, $\BB$ and $\CC$, respectively. The vector
$x_1$ can be expressed as a linear combination of the vectors
$\delta v$ ($v\in V$), which means that there is a vector $y\in\R^V$
so that $x_1=My$. Similarly, we can write $x_2=Nz$. Together with
$x=x_3$, these vectors satisfy the following system of linear
equations:
\begin{equation}\label{XYZ}
\cases{x+My+Nz=\chi_b&\cr
        M\T x=0&\cr
        N\T x=0&\cr}
\end{equation}
Multiplying the first $m$ equations by the matrix $M\T$, and using
the second equation and the fact that $M\T N=0$, we get
\begin{equation}\label{MY}
M\T My=M\T \chi_b,
\end{equation}
and similarly,
\begin{equation}\label{NZ}
N\T Nz=N\T \chi_b.
\end{equation}
Here $M\T M$ is the Laplacian of $G$ and $N\T N$ is the Laplacian of
$G^*$, and so (\ref{MY}) implies that $y=\pi_b+c\one$ for some
scalar $c$. Similarly, $z=\pi^*_b+c^*\one$ for some scalar $c'$. Thus
\[
x=\chi_b - M\T (\pi_b+c\one) - N\T (\pi^*_b+c^*\one)=\chi_b - M\T
\pi_b - N\T \pi^*_b,
\]
which is just the formula in the lemma, written in matrix form.
\end{proof}

\subsection{Nondegeneracy properties of smooth circulations}

We state and prove two key properties of smooth circulations: one,
that the projection of a basis vector to the space of smooth
circulations is non-zero, and two, that smooth circulations are
spread out essentially over the whole graph in the sense that every
connected piece of the graph where a non-zero smooth circulation
vanishes can be isolated from the rest by a small number of points.

\begin{theorem}\label{NONVANISH}
If $g>0$, then for every edge $e$, $\eta_e\not=0$.
\end{theorem}

\begin{proof}
Suppose that $\eta_e=0$. Then by Lemma \ref{SMOOTHFORM}, there are
vectors $\pi=\pi(e)\in\R^V$ and $\pi^*=\pi^*(e)\in\R^F$ such that
\begin{equation}\label{PIEQUAL}
\pi_{h(a)}-\pi_{t(a)} = \pi^*_{r(a)}-\pi^*_{l(a)}
\end{equation}
for every edge $a\not=e$, but
\begin{equation}\label{PIDIFF}
\pi_{h(e)}-\pi_{t(e)} = 1+ \pi^*_{r(e)}-\pi^*_{l(a)}.
\end{equation}

For convenience, we orient every edge so that
$\pi^*_{r(a)}\ge\pi^*_{l(a)}$.

Let $\alpha_1<\alpha_2<\dots<\alpha_k$ be the different values of
$\pi^*$, and let $\FF_i=\{F\in \FF:~\pi^*_F=\alpha_i$. Since for every
face $F$ other than $r(e)$ and $l(e)$ the value of $\pi^*_F$ is the
average of its values on the neighbors, it follows that
$\FF_1=\{l(e)\}$ and $\FF_k=\{\FF\setminus r(e)\}$.

Consider the union $S_i$ of all faces in $\FF_1\cup \dots \cup \FF_i$,
$1\le i< k$. We claim that its boundary is a single cycle containing
the edge $e$. It is trivial that $e\subseteq\partial S_i$. Suppose
that the boundary of $S_i$ contained a cycle $C$ that did not go
through $e$, which is a directed cycle with $S_i$ being on its left
hand side. But then $\pi^*_{r(a)}>\pi^*_{l(a)}$ for every edge of
$C$. By (\ref{PIEQUAL}), this means that $\pi_{h(a)}>\pi_{t(a)}$ for
every edge $a$ of $C$, which is clearly impossible.

Thus we know that $\partial S_i$ consists of a single cycle through
$e$; in other words, it consists of $e$ and a path $P_i$ connecting
$t(e)$ and $h(e)$. As above, it follows that $P_i$ is a directed path
from $h(e)$ to $t(e)$, and $\pi_v$ strictly increases as we traverse
the path.

Consider the two paths $P_{i-1}$ and $P_i$. These may have common
nodes besides their endpoints, but from the fact that $\pi$ strictly
increases along both of them, it follows that their common nodes are
in the same order on both paths. This implies that the two paths
have the following structure: there are common (directed) subpaths
$Q_0, Q_1, \dots, Q_r$ (possibly consisting of just a single node),
and cycles $C_1, C_2, \dots, C_r$ so that $C_j$ is the union of two
directed paths $R_j$ and $R'_j$, connecting the endpoint of
$Q_{i-1}$ to the starting point of $Q_i$, with $R_j\subseteq
P_{i-1}$ and $R_j'\subseteq P_i$.

The union of all faces $F\in \FF_i$ is a (not necessarily connected)
surface $S_i'$ whose boundary is the union of the cycles $C_j$.

\medskip

\noindent{\bf Claim.} {\it $S_i'$ consists of $r$ discs $D_1,\dots
D_r$, where the boundary of $D_j$ is the cycle $C_j$.}
\medskip

Let $G_i'$ be the subgraph of $G$ contained in $S_i'$ (including the
nodes on the boundary, but not the edges). Let $H$ be a connected
component of $G'$. Since every edge of $G'$ has
$\pi_{h(a)}=\pi_{t(a)}$ by (\ref{PIEQUAL}), all nodes of $H$ have the
same $\pi$ value. Since $\pi$ is strictly increasing along $P_{i-1}$
as well as along $P_i$, it follows that $H$ can be attached to at
most one node on each of these paths. By the 3-connectivity of $G$,
it follows that $H$ is a single edge connecting some $u\in
V(P_{i-1})$ and $v\in V(P_i)$. From $\pi_u=\pi_v$ it also follows
that $H$ connects a node on an $R_j$ to a node on the corresponding
$R_j'$, and the endpoints of these edges are in the same order on
both paths, which gives an ordering of the edges of $G'$.

Now every face in $\FF_i$ is a disc attached along a cycle in
$P_{i-1}\cup P_i\cup G'$, so that the edges in $G'$ have two faces
attached, while the edges of $P_{i-1}\cup P_i$ have only one. The
only way to this is to attach each face in $\FF_i$ to

\smallskip

--- two consecutive edges of $G'$ and to the to subpaths of $G$ between
their endpoints, or

\smallskip

--- to the first or last edge of $G'$ connecting two points on a cycle
$C_j$ and to the corresponding arc of $C_j$.

\smallskip

Hence the Claim follows immediately.

Now we see that $S_i$ is obtained from $S_{i-1}$ by attaching the
disjoint disks $D_1,\dots,D_r$ along single arcs $R_1,\dots,R_r$.
Since $S_1$ is a disc, it follows by induction that $S_i$ is a disk,
and in particular, $S_{k-1}$ is a disc. Since $S$ is obtained from
$S_{k-1}$ by gluing in the last face (a disc) along the cycle
$\partial S_{k-1}$, it follows that $S$ is a sphere, and so $g=0$.
\end{proof}

\begin{theorem}\label{NONDEG}
Let $G$ be a graph embedded in an orientable surface $S$ of genus $g$
so that all faces are discs. Let $h$ be a non-zero smooth circulation
on $G$ and let $G'$ be the subgraph of $G$ on which $h$ does not
vanish. Suppose that $h$ vanishes on all edges incident with a
connected subgraph $U$ of $G$. Then $U$ can be separated from $G'$ by
at most $16g$ points.
\end{theorem}

The assumption that the connectivity between $U$ and the rest of the
graph must be linear in $g$ is ``essentially'' sharp in the
following sense. Suppose that $G$ has a cutset of $2g-1$ or fewer
edges. Since the dimension of the space of smooth circulations is
$2g$, there will be a non-zero smooth circulation $h$ vanishing on
these edges. If this circulation is non-zero on (say) the left hand
side of this cut, then we can replace all

Before the proof of this theorem we need a simple lemma about maps.

\begin{lemma}\label{REVERSE}
Let $G$ be any digraph embedded on an orientable surface $S$ of genus
$g$. Assume that $G$ has no sources and sinks. For every face $F$,
let $a_F$ denote the number of nodes $v$ on the boundary of $F$ for
which the two edges on the boundary of $F$ incident with $v$ are
directed to $v$. Then
\[
\sum_F (a_F-1) \le 2g-2.
\]
\end{lemma}

\begin{proof}
Let $n$, $m$ and $f$ denote the number of nodes, edges and faces.
Clearly $a_F$ is also the number the nodes on the boundary of $F$
with both edges oriented out, and so $\sum_F 2a_F$ counts the number
of ``corners'' with both edges oriented in or both edges oriented
out. Since there are no sources or sinks, at every node there are at
least 2 corners with one edge in and one out, and hence
\[
\sum_F 2a_F \le \sum_v (d_v-2) = 2m-2n= 2f+4g-4.
\]
Rearranging and dividing by 2, we get the inequality of the lemma.
\end{proof}

\begin{corollary}\label{ONEFACE}
If there is no face whose boundary is a directed cycle, then $a_F\le
2g-1$ for every face.
\end{corollary}

Now we are ready to prove Theorem \ref{NONDEG}.

\begin{proof}
By re-orienting edges of $G'$ we may assume that $h>0$ on the edges
$G'$. Let us shrink every connected component of $G\setminus V(G')$
to a single point, to get a digraph $G''$ embedded in $S$. Note that
every face of $G''$ is a disc and it contains at least one edge of
$G'$, and so it cannot be bounded by a directed cycle.

Let $w$ be the node to which $U$ is contracted, let $F'$ be the face
of $G'$ containing $w$. We think of $F$ as a surface with a boundary
$\partial F$ consisting of $k$ disjoint Jordan curves
$C_1,\dots,C_k$, glued to some (not necessarily disjoint) cycles of
the graph $G'$.

Let $R$ be a face of $G''$ contained in $F'$, bounded by a cycle
$\partial R$. Then $\partial F\cap \partial R$ consists of one or
more arcs; let $b_R$ denote the number of these arcs.

Let $r_i$ denote the number of faces $R$ with $b_R=i$, and let $r_i'$
be the number of faces among these which touch $w$.

First we estimate $r_1$. Along the common arc of $\partial F$ and
$\partial R$, the direction of the edges must change at least once:
else, $h$ would add up to a non-zero value along the face boundary
$\partial R$, contradicting the assumption that $h$ is smooth. Hence
by Lemma \ref{ONEFACE}, $r_1\le 2a_{F'}\le 4g-2$.

Consider a face $R$ with $b_R>1$ touching $w$. Select a node $v_R$ in
the interior of $R$, and connect it inside $R$ by disjoint arcs to
$w$ and to one point on each of the arcs of $\partial F\cap
\partial R$. Call these arcs {\it red}.

Cut away all of the surface outside $F'$, and contract each $C_i$ to
a single point. Also contract $W$ to a single point. This way we get
a surface $S^*$ with genus $g^*$. The red arcs form a graph $G^*$
embedded in $S^*$. It is clear that $G^*$ is bipartite.

Next, we estimate the number $f_2^*$ digonal faces of $G^*$. Indeed,
a digonal face must be formed by two red arcs originally connecting a
node $v_R$ to two points $a$ and $a'$ on different common arcs $A$
and $A'$ of $\partial R$ and a $C_i$, and by a disc $R'$ bounded by
these two red arcs and an arc $B$ of $C_i$ connecting $a$ and $a'$.
Since $A$ and $A'$ are distinct arcs, there must be a face of $G$
inside $R'$ attached to at least one edge of $C_i$. In fact, there
must be such a face $R_0$ that is attached to $C_i$ along a single
arc. Clearly, $R_0$ is disjoint from $W$ and so it is counted in
$r_1-r_1'$. It is also clear that different digonal faces of $G^*$
correspond to different faces counted in $r_1-r_1'$. This proves that
$f_2^*\le r_1-r_1'$.

Now use Euler's formula. The number of nodes of $G^*$ is $1+k+k'$,
where the 1 accounts for $w$, $k$ is the number of components $C_i$,
and $k"=\sum_{i\ge 2} r_i'$ is the number of faces of $G$ meeting
$W$. The number of edges is $m^*=\sum_{i\ge 2} (i+1)r_i'$. So the
number of faces is
\[
f^*=\sum_{i\ge 2} (i+1)r_i'-\left(1+k+\sum_{i\ge 2}
r_i'\right)+2g^*-2.
\]
There are $f_2^*$ digons and $f^*-f_2^*$ faces that are all at least
4-gons, and hence
\[
2m^* \ge 2f_2^* + 4(f*-f_2^*)= 4f^*-2f_2^* \ge 4f^*-2(r_1-r_1')\ge
4f^*-4g-2+2r_1',
\]
and so
\[
\sum_{i\ge 2} (i+1)r_i' = m^*\ge 2 f^*-4g-2+r_1'=r_1'+\sum_{i\ge 2}
2ir_i'-2k-2-4g^*+4-4g+2
\]
whence
\[
r_1'+\sum_{i\ge 2} (i-1)r_i'\le 2k+4g^*+4g-4.
\]

Let $N$ denote the number of neighbors of $w$ on $\partial F$. Each
neighbor of $w$ is incident with at least two faces of $G''$
containing $w$ and an arc of $\partial F$, and a face containing $i$
arcs of $\partial F$ is counted at most $2i$ times this way; hence
\[
2N\le \sum_{i\ge 1}2ir_i' \le 4\left(r_1'+\sum_{i\ge 2}
(i-1)r_i'\right) \le 8k+16g+16g^*-16.
\]

To conclude the proof, it suffices to show that
\begin{equation}\label{GBOUND}
2g*+k\le 2g
\end{equation}
unless $g=g^*=0$ and $k=1$. To this end, notice that no cycle $C_i$
bounds a disk not containing $w$. Indeed, the restriction of $h$ to
the connected component of $G'$ contained in this disk would be a
smooth circulation on a planar map, contradicting Lemma
\ref{SMOOTHDIM}.

If $k=1$ and $C_1$ bounds a disc containing $w$, then $g^*=0$ and
(\ref{GBOUND}) is trivial. In every other case, none of the cycles
$C_i$ is null-homologous.

Let us cut the surface $S$ along each $C_i$ one by one, and in each
case, glue discs on both copies of $C_i$ obtained by cutting. If this
operation separates the surface, then we keep the part containing $w$
only. This way we obtain the surface $S^*$.

Every cut reduces the genus of the surface, except when it separates
the surface and the part thrown away is a sphere. Since each such
part must have been attached along at least two cycles $C_i$, at
least half of the cuts reduce the genus, showing that $g^*\le g-k/2$.
This proves (\ref{GBOUND}) and completes the proof of Theorem
\ref{NONDEG}.
\end{proof}

\section{Proof of Theorem \ref{OBSERVE}}

\subsection{A lemma about convergence}

We prove a simple lemma about the convergence of a simple randomized
iterative projection process.

Let $a_1,\dots,a_k\in\R^n$, and let $A=(1/k)\sum_{i=1}^k a_ia_i\T$.
Let $\lambda$ be the smallest positive eigenvalue of $A$. Let $L$ be
the linear subspace generated by $a_1,\dots,a_k$, and $L^\perp$, its
orthogonal complement.

For $x\in \R^n$, define a Markov chain $X^0,X^1,\dots\in\R^n$ as
follows: start with $X^0=x$. Given $X^t$, choose a vector $a_i$
(uniformly and randomly), and let
\[
X^{t+1} = X^t-{a_i\T X^t\over a_i\T a_i}a_i.
\]
Let $y$ and $z$ be the orthogonal projections of $x$ onto $L^\perp$
and $L$, respectively. Then $X^t\to y$, and in fact the following
lemma describes the rate of convergence:

\begin{lemma}\label{CONVERGE}
For every $x\in\R^n$,
\[
\E\left(|X^t-y|^2\right)\le (1-\lambda)^t|z|^2.
\]
\end{lemma}

\begin{proof}
First we consider the case $t=1$:
\[
|X^1-y|^2 = \left|z - {a_i\T z\over a_i\T a_i}a_i\right|^2 = |z|^2 -
{(a_i\T z)^2\over |a_i|^2},
\]
and hence
\[
\E\left(|X^t-y|^2\right)= |z|^2-z\T\E\left({1\over a_i\T
a_i}|a_i|^2\right)z
\]
\[
= |z|^2-z\T\left({1\over k}\sum_{i=1}^k{1\over
|a_i|^2}a_ia_i\T\right)z = |z|^2-z\T Az.
\]
Since $z\in L$ is in the range of $A$, we have $z\T Az\ge \lambda
|z|^2$, which proves the assertion.

The general case follows by induction.
\end{proof}

\subsection{Setup for the proof}

Let $a_i\in\R^E$ be defined for every node $i$ by
\[
(a_i)_e=\cases{1, & if $i$ is the head $e$,\cr
               -1, & if $i$ is the tail of $e$, \cr
               0, & otherwise.\cr}
\]
Furthermore, let $b_F\in\R^E$ be defined for every face $F$ by
\[
(b_F)_e=\cases{1, & if $e$ is an edge of $\partial F$ oriented
                    clockwise,\cr
               -1, & if $e$ is an edge of $\partial F$
                    oriented counterclockwise,\cr
               0, & otherwise.\cr}
\]
We consider the matrices
\[
A={1\over n}\sum_{i\in V} {1\over |a_i|^2} a_ia_i\T,
\]
and
\[
B={1\over f}\sum_{F\in\FF} {1\over |b_F|^2} b_Fb_F\T.
\]
Let $\lambda_1$ and $\lambda_2$ be the smallest positive eigenvalue
of $A$ and $B$, respectively, and let
$\mu=\min\{\lambda_1,\lambda_2\}$.

Let $T_t:~\R^E\to\R^E$ denote the (random) linear mapping that (a)
and (b) generate in step $t$. Note that subspace $\CC$ is invariant
under $T_t$. Let $W(t,s)=T_{s}T_{s-1}\dots T_{t+1}$,
\[
u(t)=\cases{e_j, & if in step $t$ edge $j$ was excited,\cr
            0,& otherwise.\cr}
\]
and $u(t,s)=W(t,s)u(t)$. Hence
\[
x(s)=\sum_{t=0}^s u(t,s).
\]

Let $u_1(t,s)$, $u_2(t,s)$ and $u_3(t,s)$ denote the orthogonal
projections of $u(t,s)$ to $\AA$, $\BB$ and $\CC$, respectively, and
let $u_i(t)=u_i(t,t)$. This notation is somewhat redundant, since
$\CC$ is invariant under $T_t$, and hence $u_3(t,s)=u_3(t)$ for every
$s\ge t$. The ``error part'' is $w(t,s)=u_1(t,s)+u_2(t,s)$. Thus we
get the ``smooth part''
\[
x'(s)=\sum_{t=0}^s u_3(t),
\]
and the ``error part''
\begin{equation}\label{XW}
x''(s)=\sum_{t=0}^s w(t,s).
\end{equation}
We need to bound $x'$ from below and $x''$ from above.

\subsection{Bounding smooth circulations}

The vector $u_3(t)$ is a smooth circulation. Suppose that it is not
the 0 circulation. Then by Theorem \ref{NONDEG}, we know that its
restriction $v(t)$ to $E_0$ is not the 0 circulation. We need a lower
bound on $|v(t)|$. Recall that $f = |\FF|$; then we have

\begin{lemma}\label{SMOOTHBOUND}
\[
|v(t)|\ge n^{-n}f^{-f}.
\]
\end{lemma}

\begin{proof}
$x$ is a rational vector whose denominator is a divisor of $\det(M\T
M+J)\det(N\T N+J)$. By Hadamard's inequality, the denominator of $x$
is at most $n^nf^f$.

Now $v(t)$ is a restriction of $x$ which is nonzero, and so at least
one coordinate of $v(t)$ is a non-zero rational number with
denominator at most $n^nf^f$. This proves (\ref{SMOOTHBOUND}).
\end{proof}

\subsection{Bounding the error}

We prove a bound on the ``error term''. More exactly, we fix an
integer $a>0$, and split the error into ``old errors'' and ``new
errors'':
\[
x''(s)=\sum_{t=0}^{s-a} w(t,s) + \sum_{t=s-a+1}^{s}
w(t,s)=X_1(s,a)+X_2(s,a).
\]
First we estimate the expectation of $|X_1(s)|^2$. We claim that
\begin{equation}\label{EX1}
\E(|X_1(s,a)|^2)\le 5{p\over\mu} (1-\mu)^a.
\end{equation}

Let us fix the ``excitations'' $u(0), u(1), \dots$, and let $\E_u$
denote expectation conditional on these. Lemma \ref{CONVERGE}
implies that
\[
\E_u(|u_1(t,s)|^2)\le (1-\lambda_1)^{s-t}|u_1(t)|^2,
\]
and
\[
\E_u(|u_2(t,s)|^2)\le (1-\lambda_2)^{s-t}|u_2(t)|^2,
\]
and so
\begin{eqnarray}\label{WBOUND}
\E_u(|w(t,s)|^2)&=& \E_u(|u_1(t,s)|^2+|u_2(t,s)|^2)\\
&\le&(1-\mu)^{s-t}(|u_1(t)|^2+|u_2(t)|^2)\nonumber \\
&\le& (1-\mu)^{s-t}|u(t)|^2.
\end{eqnarray}
From the definition of $\Sigma_1$ we have
\[
\E_u(|\Sigma_1(s)|^2)=\E_u\left(\left|\sum_{t=0}^{s-a}  w(t,s)
\right|^2\right)
\]
\[
=\E_u\left(\sum_{t=0}^{s-a}\sum_{t'=0}^{s-a} w(t,s)\T w(t',s)\right)=
\sum_{t=0}^{s-a}\sum_{t'=0}^{s-a} \E_u(w(t,s)\T w(t',s))
\]
\[
\le \sum_{t=0}^{s-a}\sum_{t'=0}^{s-a}
\E_u(|w(t,s)|^2)^{1/2}\E_u(|w(t',s)|^2)^{1/2}= \left(\sum_{t=0}^{s-a}
\E_u(|w(t,s)|^2)^{1/2}\right)^2.
\]
Using (\ref{WBOUND}), this gives
\begin{equation}\label{XBOUND}
\E_u(|\Sigma_1(s)|^2)\le \left(\sum_{t=0}^{s-a}
(1-\mu)^{(s-t)/2}|u(t)| \right)^2.
\end{equation}
Now we take expectation over the sequence $u(t)$. Since the $|u(t)|$
are independent 0-1 valued variables with mean $p$, the expectation
of the right hand side is easy to estimate:
\[
\E\left(\left(\sum_{t=0}^{s-a}
(1-\mu)^{(s-t)/2}|u(t)|\right)^2\right)
\]
\[
= \sum_{t=0}^{s-a} \sum_{t'=0}^{s-a}
(1-\mu)^{(s-t)/2}(1-\mu)^{(s-t')/2} \E(|u(t)||u(t')|).
\]
Here
\[
\E(|u(t)||u(t')|)=\cases{p,& if $t=t'$,\cr
                         p^2,& otherwise.\cr}
\]
Thus we can write the sum above as
\[
p^2\sum_{t=0}^{s-a} \sum_{t'=0}^{s-a}
(1-\mu)^{(s-t)/2}(1-\mu)^{(s-t')/2} + (p-p^2)\sum_{t=0}^{s-a}
(1-\mu)^{s-t},
\]
\[
=p^2(1-\mu)^a\left(1-(1-\mu)^{(s-a)/2}\over 1 -
(1-\mu)^{1/2}\right)^2 + (p-p^2)(1-\mu)^a {1-(1-\mu)^a\over
1-(1-\mu)}
\]
\[
\le (1-\mu)^a {4p^2\over \mu^2} + {p-p^2 \over \mu} < 5{p\over\mu}
(1-\mu)^a.
\]
This proves (\ref{EX1}).

Second, we consider $X_2$. If no excitation event occurs between
times $s-a+1$ and $s$, then $\Sigma_2(s)=0$. Hence
\begin{equation}\label{PX2}
\Pr(|\Sigma_2(s)|>0) \le 1-(1-p)^a <ap.
\end{equation}

From (\ref{EX1}) and (\ref{PX2}) we get that
\[
\Pr(|x''(s)|>\delta) \le \Pr(|X_1(s)|\ge \delta) + \Pr(|X_2(s)|>0)
<   {5p\over\mu\delta^2} (1-\mu)^a + ap.
\]
The choice $a={2\over\mu}\ln{1\over\delta}$ gives the best bound (up
to a constant):
\begin{equation}\label{ERRORBOUND}
\Pr(|x''(s)|> \delta) < {10p\over\mu}\ln{1\over\delta}.
\end{equation}

\subsection{Completing the proof}

To simplify notation, put $z_i=z(t_i)$, $z'_i=z'(t_i)$ and
$z''_i=z''(t_i)$. It follows from the choice of the integers $t_i$
that the vectors $z_1, z_2, \dots z_k$ are linearly independent, and
so $\dim(\LL(i))=i$. Furthermore, by the selection of these vectors,
the Gram-Schmidt orthogonalization $z_1^*=z_1, z_2^*, \dots, z_k^*$
consists of vectors of length at least $\eps$.

Let $\delta=(\eps/4)(1+1/\eps)^{-2g_0}$. Since $|z''(t)|\le|x''(t)|$,
the probability that $|z''(t)|>\delta$ is less than
${10p\over\mu}\ln{1\over\delta}$ by (\ref{ERRORBOUND}), so the
probability that any of the $t_i$ have $|z''(t_i)|>\delta$ is less
than
\[
g_0{10p\over\mu}\ln{1\over\delta}<{1\over 9}.
\]
So with probability at least $.99$, we have $|z''_i|\le \delta$ for
every $i$. Let us assume that this occurs.

\smallskip

\noindent{\bf Claim 1.} Suppose that for some real numbers
$\alpha_1,\dots,\alpha_k$, we have
\begin{equation}\label{ASMALL}
\left|\sum_{i=1}^k \alpha_i z_i\right|\le 1.
\end{equation}
Then
\[
|\alpha_i|\le {1\over\eps}\left(1+{1\over \eps}\right)^{k-i}.
\]

\smallskip

Indeed, (\ref{ASMALL}) implies that for any $1\le j\le k$,
\[
\left|\sum_{i=1}^k\alpha_i z_i\T z_j^*\right|\le |z_j^*|.
\]
Since $z_i\T z_j^*=0$ if $i<j$ and $z_j\T z_j^*=|z_j^*|^2$, it
follows that
\[
|\alpha_j||z_j^*|^2 \le |z_j^*| - \sum_{i=j}^k |\alpha_iz_i\T
z_j^*|\le |z_j^*| + \sum_{i=j}^k |\alpha_i||z_j^*|.
\]
Dividing by $|z_j^*|^2$ and using that $|z_i|\le 1$, we get that
\[
|\alpha_j|\le {1\over |z_j^*|}\sum_{i=j}^k |\alpha_i|\le {1\over
\eps}\sum_{i=j}^k |\alpha_i|.
\]
Hence the claim follows by induction on $k-i$.

\smallskip

\noindent{\bf Claim 1.} The vectors $z'_1,\dots,z'_{2g}$ are linearly
independent.

\smallskip

Indeed, assume that there is a linear relation $\sum_i\alpha_i
z'_i=0$ where not all the $\alpha_i$ are 0. We can write this as
$\sum_i\alpha_i z_i=\sum_i \alpha_i z''_i$. Since the $z_i$ are
linearly independent, the left hand side is non-zero, and so we may
assume it has norm 1. But then Claim 1 implies that
$|\alpha_i|\le(1/\eps)(1+1/\eps)^{k-i}$, and so using that
$|z_i|\le\delta$,
\[
\left|\sum_i \alpha_i z''_i\right|\le \delta\sum_i |\alpha_i|\le
\delta\left(1+{1\over\eps}\right)^k < 1,
\]
a contradiction.

\noindent{\bf Claim 3.}  $d(z(t),\LL(k))<|z''(t)|+{\eps\over 2}$.

\smallskip

Indeed, by Claim 2 we can write $z'(t)=\sum_{i=1}^k \alpha_i z'_i$
with some real numbers $\alpha_i$. Then
\[
\sum_{i=1}^k \alpha_iz_i= z'(t)+\sum_{i=1}^k \alpha_iz_i''.
\]
Let $R$ be the norm of the vector on the two sides. Then by Claim 1,
we have $|\alpha_i|\le R(1/\eps)(1+1/\eps)^{k-i}$. Using this, we get
\[
R=\left|z'(t)+\sum_{i=1}^k \alpha_iz_i''\right|\le 1 +
R\delta\left(1+{1\over\eps}\right)^k \le 1 + {1\over 2}R,
\]
whence $R\le 2$, and so $|\alpha_i|\le (2/\eps)(1+1/\eps)^{k-i}$.

Let $w=\sum_{i=1}^k\alpha_i z_i\in \LL(k)$. Then
\[
|w-z'(t)|=\left|\sum_{i=1}^k \alpha_iz_i''\right| \le
\delta\sum_{i=1}^k |\alpha_i| \le 2\delta
\left(1+{1\over\eps}\right)^k <\eps/2,
\]
and so
\[
d(z(t),\LL(k))\le |z(t)-w|\le |z(t)-z'(t)|-|z'(t)-w| < |z''(t)|+
{\eps\over 2},
\]
as claimed.

\smallskip

Now we are ready to complete the proof. In guessing that $k/2$ is the
genus of the surface, we can err in two directions: it may be that
$g<k/2$ or $g>k/2$. We estimate the probability of these errors
separately.

If $g<k/2$, then we moved on when we had $k=2g$, which means that at
that stage we had too many vectors $z(t)$ farther from $\LL(k)$ than
$\eps$. By Claim 3, every such $t$ must satisfy $|z''(t)|>\eps/2
>\delta$. By (\ref{ERRORBOUND}), the probability that a given $t$
has this property is less than ${10p\over\mu}\ln{1\over\delta}$, and
so using Chernoff's inequality, the probability that the number of
indices $t$ with this property is at least $T'$ is less than $1/9$.

If $g>k/2$, then we stopped too early: there was a subspace $L'$
with $\dim(L')<\dim(L)$ so that almost all the $z(t)$ were closer to
$L'$ than $n^{-10n}$. The mapping $\R^E\to \LL$ obtained by
projecting $\R^E$ to $\CC$ and then restricting it to $E_0$ is
surjective; this implies that there are $2g$ edges $e_1,\dots,e_{2g}
\in E$ so that the vectors $\mu(e_1),\dots,\mu(e_{2g})$ form a basis
in $\LL$. These vectors are rational with denominators at most
$n^nf^f$, so the determinant of this basis is at least
$n^{-m_0n}f^{-m_0f}\ge \eps$. It follows that at least one of these
vectors, say $\mu(e_1)$, is at distance at least $\eps$ from $\LL$.

Thus if in any of the time intervals $[(t-1)N+1,tN]$, the edge $e_1$
was excited, no other edge was excited, and the excitation of $e_i$
occurred in the first $N-a$ steps in this time interval, then $z(t)$
must be farther from $\LL'$ than $\eps$. The probability that this
happens for a given $t$ is $(N-a)p(1-p)^{N-1}/m$; so (using the
Chernoff bound again) the probability that this happens for fewer
than $T'$ choices of $t$ is less that $1/9$.

To sum up, the total probability of ``bad'' cases is less than $1/9$
(when $|z_i|>\delta$ for some $i$) plus $1/9$ (when $g<k/2$) plus
$1/9$ (when $g>k/2$). This proves the theorem.

\section{Concluding remarks}

\noindent{\bf 1.} We can make some cosmetic changes to the setup as
given above. One objection may be that the noisy circulator, as
defined, is not truly local, since (say) in operation (a) we have to
select a node uniformly from all nodes. The standard way of fixing
this is to attach an ``alarm clock'' to each node, edge, and face,
which wakes them up at random times according to a Poisson process
(the edge-clock is much slower than the other two).

Another objection is that the noisy circulator as constructed above
is not stationary: the total mass grows to infinity. An easy fix is
to give a second, even slower clock to each edge: when this rings,
they reset their value to 0. Another possible fix comes from the
observation that the excitations don't necessarily have to be
constants, the proof works just as well for random and symmetric
excitations. So modifying the excitation step to reset the value of
any edge to $1$ with very small probability (rather than to add $1$)
provides a stationary version (but the analysis becomes more
complicated).

Further variants, improvements and generalizations of the above
system are of interest:

\begin{itemize}
\item Can one recover the genus even if the rate at which excitations
take place is faster (ideally, independent of the number of nodes)?

\item Suppose that we only allow two values (or any other given
discrete set of values) on any edge. Can we still recover the genus?

\item Can one extend the technique to recover a non-orientable surface
by observing a random process on an embedded graph locally?
\end{itemize}

For  background regarding graphs on surfaces and Riemann surfaces
see for instance \cite{MT} \cite{Bu}, for background regarding algorithmic
applications of random processes see  \cite{H}.

\medskip

\noindent{\bf 2.} The notion of global information from local
observation can inspire many questions in different directions. We
briefly present some related examples.

\begin{example}\label{E1} Let $G$ be a finite connected graph. Start a simple
random walk on $G$. Fix a vertex $v$ in $G$. You are given the
sequence of times for which the simple random walk visits $v$, what
information can be learned about $G$? From the infinite sequence one
can reconstruct the on diagonal heat kernel and thus the spectrum of
the transition matrix of the random walk (so in the case of regular
graphs, the spectrum of the graph).
\end{example}

\begin{example}\label{E2} [Obtaining the size of the road system by
measuring the volume of local traffic]. Let $G$ be a finite connected
irreducible regular graph. From each vertex start an independent
simple random walk. Fix a vertex $v \in G$. For each time $t$ you
can observe the number of walkers occupying $v$. How much time
is needed in order to, almost surely, know the size of $G$? Assume
you are given an {\it a priori} bound $N$ on $n=|V(G)|$. Then one
can get a (poor) polynomial upper bound on the time needed along the
following lines. The mixing time of $G$ is smaller than $n^3$, so if
we observe the load on $v$ only at times of the form $kN^3$, we get
almost independent samples from a distribution which is
exponentially close to Binomial$(n, 1/n)$, the stationary
distribution for the number of walkers at $v$. For Binomial$(n,
1/n)$, about $n^4$ samples are needed to recover $n$. So after about
$N^7$ steps $n$ can be recovered. (Probably $N^4$ is the time needed
to recover $n$.)
\end{example}

There are several random processes on graphs which have been
considered before, and for which the question ``what global
information can be deduced from local observation'' is meaningful.
We mention two examples.

\begin{example}\label{E3} Consider an $n$-vertex connected graph. There
are $k \leq n$ particles labeled $1,2,...,k$. In a configuration,
there is one particle at each vertex. The interchange process
discussed briefly in \cite {A}, is a continuous-time Markov chain on
configurations. For each edge  $(i,j)$, at rate $1$ the particles at
vertex $i$ and vertex $j$ are interchanged. If only one of the two
vertices is occupied, then it jumps to the other vertex. Assume you
observe which labeled particle occupies a fixed predetermined vertex
at any time. For $k=1$, this is just Example \ref{E1}. If $k >1$, can
one recover, using these observations, further information about $G$
not contained in the spectrum? If $k = n$, is it possible to
reconstruct $G$?
\end{example}

\begin{example}\label{E4}
Another natural candidate for local observation is the heat-bath
chain (Glauber dynamics) on $k$-colorings of a graph: at any step,
we have a (legal) $k$-coloring; we select a random node $v$ and a
random color $\alpha$, and we re-color $v$ with color $\alpha$ if
this gives a legal $k$-coloring. Can we derive estimates on the
chromatic number, or maximum degree, by observing a bounded piece of
the graph?
\end{example}

\medskip

\noindent{\bf Acknowledgement.} We are grateful to Dimitris
Achlioptas, Mike Freedman, Eran Makover, Oded Schramm and Kevin
Walker for stimulating discussions on this topic.

\end{document}